\documentclass[11pt,reqno]{amsart}
\usepackage{amssymb}
\usepackage[hidelinks]{hyperref}
\hypersetup{pdfauthor={Bo Yang},pdftitle={Finite-time blowup for the 3D focusing cubic nonlinear Schr\"odinger equation}}
\numberwithin{equation}{section}
\newtheorem{theorem}{Theorem}[section]
\newtheorem{proposition}[theorem]{Proposition}
\newtheorem{lemma}[theorem]{Lemma}
\newcommand{\R}{\mathbb R}
\newcommand{\dd}{\,d}
\newcommand{\Qprod}{\|Q\|_{L^2_x}\|\nabla Q\|_{L^2_x}}
\newcommand{\Sch}{\mathcal S}
\newcommand{\HL}{\mathcal M}
\newcommand{\im}{\operatorname{Im}}
\newcommand{\re}{\operatorname{Re}}
\newcommand{\supp}{\operatorname{supp}}
\newcommand{\ip}[2]{\langle #1,#2\rangle}

\title[Finite-time blowup for the focusing cubic NLS]{Finite-time blowup for the 3D focusing cubic nonlinear Schr\"odinger equation}
\author[Bo Yang]{Bo Yang}
\address{School of Mathematics and Statistics, Huazhong University of Science and Technology, Wuhan 430074, Hubei, P. R. China; Hubei Key Laboratory of Engineering Modeling and Scientific Computing, Huazhong University of Science and Technology, Wuhan 430074, Hubei, P. R. China}
\email{boyang\_0717@hust.edu.cn}
\date{}

\begin{document}
\begin{abstract}
We prove the Holmer--Roudenko strong conjecture for the 3D focusing
cubic nonlinear Schr\"odinger equation $i\partial_t u+\Delta u+|u|^2u=0$. 
More precisely, let $M(u)$ and $E(u)$ denote the mass
and energy of a solution $u$, and let $Q$ be the ground state solution
of $-Q+\Delta Q+|Q|^2Q=0$. For any initial data $u_0\in H^1(\mathbb{R}^3)$
satisfying $M(u)E(u)<M(Q)E(Q)$ and
$\|u_0\|_{L^2_x}\|\nabla u_0\|_{L^2_x}>\|Q\|_{L^2_x}\|\nabla Q\|_{L^2_x}$,
we prove that $u$ blows up in finite time in both time directions,
without any radial symmetry or finite-variance assumption. Moreover, we
obtain a quantitative upper bound $C_\delta M(u)^2$ on the length of the
maximal interval of existence, where $C_\delta>0$ depends only on the
mass--energy gap $\delta=1-(M(u)E(u))_+/(M(Q)E(Q))$. The proof
introduces a projection--averaging virial method, combining
slope-constrained $L^2$ projection of the cumulative half-space mass
with logarithmic averaging to construct a bounded, mass-adaptive virial
multiplier.
\end{abstract}
\maketitle

\section{Introduction}

We consider the three-dimensional focusing cubic nonlinear Schr\"odinger
equation (NLS)
\begin{equation}\label{eq:nls}
 i\partial_t u+\Delta u+|u|^2u=0,\qquad
 u(x,0)=u_0(x),\qquad u_0\in H^1(\R^3),
\end{equation}
where $u=u(x,t)$ is complex-valued. Solutions to \eqref{eq:nls} in
$H^1(\R^3)$ conserve mass and energy, defined by
\[
 M(u)=\|u\|_{L^2_x}^2,\qquad
 E(u)=\frac12\|\nabla u\|_{L^2_x}^2-\frac14\|u\|_{L^4_x}^4.
\]
By conservation, we write $M(u)=M(u(t))=M(u_0)$ and
$E(u)=E(u(t))=E(u_0)$.

Equation \eqref{eq:nls} is invariant under the scaling
$u_\lambda(x,t)=\lambda u(\lambda x,\lambda^2 t)$ for $\lambda>0$, which
preserves the $\dot H^{1/2}(\R^3)$ norm; hence \eqref{eq:nls} is
mass-supercritical and energy-subcritical. Under this scaling,
\[
 M(u_\lambda)=\lambda^{-1}M(u),\qquad
 E(u_\lambda)=\lambda E(u),\qquad
 \|\nabla u_\lambda(0)\|_{L^2_x}=\lambda^{1/2}\|\nabla u_0\|_{L^2_x},
\]
so both $M(u)E(u)$ and $\|u_0\|_{L^2_x}\|\nabla u_0\|_{L^2_x}$ are
scale-invariant.

The local $H^1$ theory based on Strichartz estimates \cite{Caz,CW}
gives a unique maximal solution $u\in C(I;H^1(\R^3))$ on
$I=(-T_*,T^*)$, where $T_*,T^*\in(0,\infty]$. The solution depends
continuously on $u_0$ on compact subintervals of $I$. If $T^*<\infty$
(resp.\ $T_*<\infty$), then $\|\nabla u(t)\|_{L^2_x}\to\infty$ as
$t\uparrow T^*$ (resp.\ $t\downarrow -T_*$). If $T^*=\infty$ (resp.\
$T_*=\infty$), the solution is global forward (resp.\ backward) in time.

Consider the stationary elliptic equation
\[
 -\phi+\Delta\phi+|\phi|^2\phi=0,\qquad x\in\R^3 .
\]
It admits a smooth, positive, radial ground state $Q\in H^1(\R^3)$
that decays exponentially at infinity and has minimal mass among all
nonzero $H^1$ solutions \cite{Caz,Weinstein}, and is the unique
positive radial $H^1$ solution \cite{Coffman}. It gives the standing-wave solution
$u(x,t)=e^{it}Q(x)$ of \eqref{eq:nls}, which exists globally in time
and neither scatters nor blows up.

When $M(u)E(u)<M(Q)E(Q)$ and
$\|u_0\|_{L^2_x}\|\nabla u_0\|_{L^2_x}<\Qprod$, the solution is global
and scatters in both time directions. This was proved for radial data by
Holmer and Roudenko \cite{HRrad} and extended to general $H^1$ data by
Duyckaerts, Holmer, and Roudenko \cite{DHR}. At the exact threshold
$M(u)E(u)=M(Q)E(Q)$, the dynamics of solutions were classified in
\cite{DR10} for \eqref{eq:nls} (see also \cite{CFR} for general
mass-supercritical, energy-subcritical nonlinearities).

\subsection{Blowup above the mass--gradient threshold}
We consider initial data satisfying
\begin{equation}\label{eq:hyp}
 M(u)E(u)<M(Q)E(Q),\qquad
 \|u_0\|_{L^2_x}\|\nabla u_0\|_{L^2_x}>\Qprod .
\end{equation}
For negative energy, Glassey \cite{Glassey} proved finite-time blowup
for initial data with finite variance ($|x|u_0\in L^2(\R^3)$), and Ogawa
and Tsutsumi \cite{OT} proved blowup for radial $H^1$ data using a
localized virial identity. Both results extend to the regime
\eqref{eq:hyp}, which includes positive energies: Kuznetsov, Juul
Rasmussen, Rypdal, and Turitsyn \cite{KRRT} established the
finite-variance criterion, and Holmer and Roudenko \cite{HRrad} proved
the radial case (see also \cite{DR15} for finite-variance scattering and
blowup criteria beyond the threshold $M(u)E(u)\ge M(Q)E(Q)$). Bellazzini and Forcella \cite[Theorem~1.3]{BF} also proved finite-time
blowup in both time directions under \eqref{eq:hyp} when $u_0$ is
cylindrically symmetric about the $x_3$-axis and
$x_3u_0\in L^2(\R^3)$.

Without symmetry or variance assumptions, however, finite-time blowup
under \eqref{eq:hyp} has remained open. For arbitrary $H^1$ data
satisfying \eqref{eq:hyp}, Holmer and Roudenko \cite[Theorem~1.1]{HR}
proved a ``blowup or grow-up'' alternative: either $T^*<\infty$, or
$T^*=\infty$ and there exists a sequence $t_n\to+\infty$ such that
$\|\nabla u(t_n)\|_{L^2_x}\to\infty$, with the same alternative in
negative time. The forward-global alternative is called \emph{grow-up}.
Without symmetry or finite variance, localized virial arguments also
give blowup or grow-up results in more general settings, such as
higher-power nonlinearities \cite{DWZ} and threshold solutions
\cite{GI}. In \cite{HR}, Holmer and Roudenko proposed two conjectures
under \eqref{eq:hyp}:

\smallskip
\noindent\textbf{Weak conjecture.} If $T^*=\infty$, then
$\lim_{t\to+\infty}\|\nabla u(t)\|_{L^2_x}=+\infty$.

\smallskip
\noindent\textbf{Strong conjecture.} Every such solution blows up in
finite forward time, that is, $T^*<\infty$.

\smallskip
The purpose of this paper is to prove the strong conjecture of Holmer
and Roudenko. In fact, we establish a quantitative bound on the total
lifespan that depends only on the mass and the mass--energy gap.

\begin{theorem}\label{thm:main}
Let $u$ be the maximal $H^1(\R^3)$ solution of \eqref{eq:nls}, and
assume \eqref{eq:hyp}.
Set $\delta=1-(M(u)E(u))_+/(M(Q)E(Q))$, where $s_+:=\max\{s,0\}$. There
exists a constant $C_\delta>0$, depending only on $\delta$, such that
\begin{equation}
 T_*+T^*\le C_\delta M(u)^2.                                    \label{eq:lifespan}
\end{equation}
In particular, $T_*,T^*<\infty$, and $\|\nabla u(t)\|_{L^2_x}\to\infty$
both as $t\uparrow T^*$ and as $t\downarrow -T_*$.
\end{theorem}

Observe that \eqref{eq:lifespan} is invariant under the scaling
$u\mapsto u_\lambda$, since $T_*+T^*$ and $M(u)^2$ both scale by
$\lambda^{-2}$ while $\delta$ is unchanged.

\subsection{Main ideas: the projection--averaging virial method}
We first recall the negative virial bound forced by \eqref{eq:hyp}.
Weinstein's sharp Gagliardo--Nirenberg inequality \cite{Weinstein}
gives, for all $v\in H^1(\R^3)$,
\[
 \|v\|_{L^4_x}^4\le C_{\mathrm{GN}}\|v\|_{L^2_x}\|\nabla v\|_{L^2_x}^3,
 \qquad C_{\mathrm{GN}}=\frac{4}{3\Qprod}.
\]
Recall also the Pohozaev identity
$6M(Q)E(Q)=\|Q\|_{L^2_x}^2\|\nabla Q\|_{L^2_x}^2$ (see \cite[Section~2]{HR}).
If $\|v\|_{L^2_x}\|\nabla v\|_{L^2_x}=\Qprod$, then Weinstein's inequality
gives $M(v)E(v)\ge M(Q)E(Q)$. Under \eqref{eq:hyp}, it follows from
continuity of the flow and conservation of mass and energy that
$\|u(t)\|_{L^2_x}\|\nabla u(t)\|_{L^2_x}>\Qprod$ for all $t\in I$. By
the Pohozaev identity, this gives the first inequality in
\eqref{eq:trap} below. Define the virial functional by
\[
 K(v)=\|\nabla v\|_{L^2_x}^2-\frac34\|v\|_{L^4_x}^4 .
\]
With $\delta=1-(M(u)E(u))_+/(M(Q)E(Q))\in(0,1]$ as in
Theorem~\ref{thm:main}, we have $1-\delta\ge0$ and
$M(u)E(u)\le(1-\delta)M(Q)E(Q)$. Using the lower bound on
$\|\nabla u(t)\|_{L^2_x}^2$ together with
$2K(u(t))=6E(u)-\|\nabla u(t)\|_{L^2_x}^2$, we obtain
\[
 6E(u)\le(1-\delta)\frac{6M(Q)E(Q)}{M(u)}
      \le(1-\delta)\|\nabla u(t)\|_{L^2_x}^2,
\]
and hence, for all $t\in I$,
\begin{equation}
 \|\nabla u(t)\|_{L^2_x}^2\ge\frac{6M(Q)E(Q)}{M(u)},\qquad
 2K(u(t))\le-\delta\|\nabla u(t)\|_{L^2_x}^2 .                \label{eq:trap}
\end{equation}

For finite-variance data, the virial quantity
$\im\int_{\R^3}\bar u\,x\cdot\nabla u\dd x$ is well-defined and has
time derivative $2K(u(t))$, which gives finite-time blowup by
\eqref{eq:trap}. For general $H^1$ data, this integral need not be
finite. Cutting off $x$ at a fixed radius $R$ gives a finite quantity,
but the bound on its time derivative contains an $O(R^{-2}M(u))$ error
and a tail proportional to $\int_{|x|>R}|u(x,t)|^4\dd x$
\cite[Section~3]{HR}. For radial data, the radial Gagliardo--Nirenberg
inequality controls this tail uniformly in $t$ for large $R$. For
general data, however, translating a profile preserves \eqref{eq:hyp}
and can move almost all of its $L^4$ norm outside any fixed ball.
Moreover, along a grow-up sequence $t_n\to+\infty$ with
$\|\nabla u(t_n)\|_{L^2_x}\to\infty$, energy conservation gives
\[
 \frac{\|u(t_n)\|_{L^4_x}^4}{\|\nabla u(t_n)\|_{L^2_x}^2}
 =2-\frac{4E(u)}{\|\nabla u(t_n)\|_{L^2_x}^2}\longrightarrow 2 .
\]
That is, $\|u(t_n)\|_{L^4_x}^4$ is comparable to
$\|\nabla u(t_n)\|_{L^2_x}^2$. Therefore, if a cutoff misses even a
small fixed fraction of the $L^4$ term, the resulting error can be as
large as the negative term $-\delta\|\nabla u(t)\|_{L^2_x}^2$ in
\eqref{eq:trap} and cannot be absorbed.

This shows that a fixed spatial cutoff cannot track non-radial mass
when the solution translates or splits in space. To avoid this problem,
we construct a bounded virial multiplier that adapts to the mass
distribution of the solution; we call this approach the
\emph{projection--averaging virial method}.

The first idea is to control the $L^4$ loss by planar mass densities,
independent of where the mass lies in space. Assuming first that $u$ is
a Schwartz solution, fix a coordinate $z=x_k$, write $y\in\R^2$ for the
other two coordinates, and let $m(z,t)=\int_{\R^2}|u(z,y,t)|^2\dd y$.
On each plane $x_k=z$, the two-dimensional Gagliardo--Nirenberg
inequality gives
\[
 \int_{\R^2}|u(z,y,t)|^4\dd y
 \le C\,m(z,t)\int_{\R^2}|\nabla_y u(z,y,t)|^2\dd y .
\]
Thus, once we fix a small level $\theta_1>0$, the $L^4$ contribution
from all planes with $m<\theta_1$ is bounded by
$C\theta_1\|\nabla_y u(t)\|_{L^2_x}^2$, regardless of their location.
For the directional virial quantity $J_k(t)=\int_{\R}Gp\dd z$, where
$p=\im\int_{\R^2}\bar u\,\partial_z u\dd y$, it therefore suffices to
construct a bounded multiplier $G(z,t)$ whose slope $h=\partial_z G$
lies in $[0,1]$ and satisfies $h=1$ on $\{m\ge\theta_1\}$.

The main difficulty with a time-dependent multiplier $G(z,t)$ is that
differentiating $J_k(t)$ produces an extra term $\int_{\R}G_t p\dd z$.
For a general time-dependent multiplier, this term has no definite sign and
is hard to control. Our key idea is to construct $G$ directly from the
cumulative half-space mass $F(z,t)=\int_{-\infty}^z m(s,t)\dd s$. By
the continuity equation, $F_t=-2p$, and hence
\[
 \int_{\R}G_t(z,t)p(z,t)\dd z
 =-\frac12\int_{\R}G_t(z,t)F_t(z,t)\dd z .
\]
To make this term nonpositive, it suffices that the map $F\mapsto G$ be
monotone in $L^2_z$. A natural monotone operator with a slope
constraint is the $L^2_z$ projection onto a closed convex set: for each
level $\theta>0$, we project $F/\theta$ onto the set of functions with
distributional slope at most $1$, obtaining a unique projection
$g_\theta$. This variational step achieves three things at once:
$g_\theta$ is bounded (so no finite-variance assumption is needed), its
slope equals $1$ wherever $m>\theta$ (via the least concave majorant of
a primitive, see \eqref{eq:slope}), and comparing the variational
inequalities \eqref{eq:var} at two times gives the needed monotonicity.

For a single level $\theta$, the slope $g_\theta'$ may jump where the
two branches of \eqref{eq:slope} meet, which prevents controlling the
derivative error $\int_{\R}\partial_z h\,\partial_z m\dd z$ in the
virial identity. To smooth the slope while keeping monotonicity, we
average $g_\theta$ logarithmically over $[\theta_0,\theta_1]$ with
respect to $\dd\theta/(L\theta)$, where $0<\theta_0<\theta_1$ and
$L=\log(\theta_1/\theta_0)>1$. The resulting multiplier $G$ stays
bounded, preserves $\int_{\R}G_t p\dd z\le0$ and $h=1$ on
$\{m\ge\theta_1\}$, and has a Lipschitz slope $h=\partial_z G$ whose
derivative error is bounded by $O(L^{-1})\|\partial_k u(t)\|_{L^2_x}^2$
(via a Fisher-type estimate in Lemma~\ref{lem:fisher} and
Lemma~\ref{lem:log}). Importantly, $\theta_1$ and
$L$ can be chosen independently: we first choose $\theta_1$ small to
control the $L^4$ loss on $\{m<\theta_1\}$, and then choose $L$ large
to absorb the derivative error.

Summing $J_k$ over $k=1,2,3$ gives a locally absolutely continuous
function $J(t)$ satisfying
\[
 J'(t)\le-\frac12\delta\|\nabla u(t)\|_{L^2_x}^2,\qquad
 |J(t)|\le c_\delta M(u)^{3/2}\|\nabla u(t)\|_{L^2_x}
\]
(see Proposition~\ref{prop:virial}). Combining these bounds with
\eqref{eq:trap} yields $-J'(t)\ge c_1+c_2 J(t)^2$ for constants
$c_1,c_2>0$, so $\arctan(cJ(t))$ (with $c=\sqrt{c_2/c_1}$) decreases at
least at a positive constant rate. Since
$\arctan$ is bounded, this gives the lifespan bound \eqref{eq:lifespan}
for Schwartz solutions, which then passes to general $H^1$ data by
approximation.

To the best of the author's knowledge, this projection--averaging virial method is new in
the blowup theory of nonlinear Schr\"odinger equations. This approach
combines the slope-constrained $L^2$ projection of the cumulative mass
with logarithmic averaging, yielding a mass-adaptive multiplier whose
time-derivative term is nonpositive. Unlike proofs using symmetry,
finite variance, or concentration--compactness, our argument is direct
and quantitative, giving the explicit lifespan bound
\eqref{eq:lifespan}. The method relies only on the continuity equation,
lower-dimensional Gagliardo--Nirenberg bounds on slices, and a negative
virial functional. It should therefore apply to other
mass-supercritical and energy-subcritical focusing dispersive models without radial symmetry or finite-variance assumption.

The rest of the paper is organized as follows. Section~\ref{sec:multiplier}
carries out the one-dimensional construction of $G$ via the
slope-constrained projection (Subsection~\ref{subsec:proj}) and
logarithmic averaging (Subsection~\ref{subsec:log}).
Section~\ref{sec:virial} applies $G$ to Schwartz solutions of
\eqref{eq:nls}, proves the time-derivative bound
$\int_{\R}G_t p\dd z\le0$ (Subsection~\ref{subsec:time}), and derives
the directional virial identity to prove Proposition~\ref{prop:virial}
(Subsection~\ref{subsec:jk}). Finally, Section~\ref{sec:lifespan}
proves Theorem~\ref{thm:main}, first for Schwartz solutions and then
for general $H^1(\R^3)$ data by approximation.

\section{Construction of the one-dimensional multiplier}\label{sec:multiplier}

In this section we work in a purely one-dimensional setting and
construct the virial multiplier associated with a given nonnegative
Schwartz mass density $m\in\Sch(\R)$ having total mass
$M=\int_{\R}m(z)\dd z>0$ and $\sqrt m\in H^1(\R)$. We first construct
the slope-constrained projection $g_\theta$ at a single level $\theta>0$
and identify its slope, and then average $g_\theta$ logarithmically in
$\theta$ to obtain a Lipschitz slope with a small derivative error.

\subsection{The slope-constrained projection}\label{subsec:proj}
Define the cumulative mass of $m$ by
\[
 F(z)=\int_{-\infty}^z m(s)\dd s .
\]
For any $\theta>0$, the rescaled cumulative mass $F/\theta$ tends to $0$
as $z\to-\infty$ and to $M/\theta$ as $z\to+\infty$. Choose a smooth
function $s_\theta\in C^\infty(\R)$ that equals $0$ for all
sufficiently small $z$ and $M/\theta$ for all sufficiently large $z$.
By the rapid decay of $m$, we have $F/\theta-s_\theta\in L^2(\R)$. We
therefore define the real affine space
$s_\theta+L^2(\R)=\{s_\theta+w:\ w\in L^2(\R;\R)\}$, equipped with the
$L^2_z$ distance $(v_1,v_2)\mapsto\|v_1-v_2\|_{L^2_z}$. This affine space
is independent of the choice of $s_\theta$, since any two choices differ
by a function in $C_c^\infty(\R)$.
Consider the subset
\[
 \mathcal C_\theta=\{v\in s_\theta+L^2(\R):\ v'\le1
                        \text{ in the sense of distributions}\}.
\]
The constraint $v'\le1$ means that
$-\int_{\R}v\varphi'\dd z\le\int_{\R}\varphi\dd z$ for every
nonnegative $\varphi\in C_c^\infty(\R)$. Hence $\mathcal C_\theta$ is a
closed convex subset of $s_\theta+L^2(\R)$. It is nonempty, since
\[
 v_0(z)=\min\{\max\{z,0\},M/\theta\}
\]
satisfies $v_0'\le1$ in the sense of distributions, and
$v_0-s_\theta$ is bounded and compactly supported.
By the projection theorem in Hilbert space, there is
a unique projection $g_\theta\in\mathcal C_\theta$ minimizing
$\|v-F/\theta\|_{L^2_z}$ over $\mathcal C_\theta$. Moreover, $g_\theta$
is uniquely characterized within $\mathcal C_\theta$ by the variational
inequality
\begin{equation}
 \ip{F/\theta-g_\theta}{v-g_\theta}_{L^2_z}\le0
 \qquad (v\in\mathcal C_\theta).                              \label{eq:var}
\end{equation}

To locate the intervals where $g_\theta$ has slope one, we define, for
$\sigma\in\{-1,1\}$ and $r>0$, the one-sided weighted averages and their
supremum
\begin{equation}
 U_{\sigma,r}(z)=\frac{2}{r^2}\int_0^r(r-s)m(z+\sigma s)\dd s,
 \qquad m_{\triangle}(z)=\sup_{\sigma,r} U_{\sigma,r}(z).                 \label{eq:triangular}
\end{equation}
The function $m_{\triangle}$ satisfies the following elementary
properties:
\begin{enumerate}
 \item $m_{\triangle}\ge m$ on $\R$: since the weight $2(r-s)/r^2$ has
 total integral $1$ on $[0,r]$, letting $r\downarrow0$ gives
 $U_{\sigma,r}(z)\to m(z)$ and hence $m_{\triangle}\ge m$.
 \item $m_{\triangle}>0$ on $\R$: since $m\ge0$ and $M>0$, for every
 $z\in\R$ we can choose $\sigma\in\{-1,1\}$ and $r>0$ such that
 $\supp m$ intersects the open interval between $z$ and $z+\sigma r$,
 giving $U_{\sigma,r}(z)>0$.
 \item $m_{\triangle}(z)\to0$ as $|z|\to\infty$: for any $R>0$, we have
 $U_{\sigma,r}(z)\le 2M/r\le 2M/R$ for $r\ge R$ and
 $U_{\sigma,r}(z)\le\sup_{|w-z|\le R}m(w)$ for $r<R$; letting
 $|z|\to\infty$ followed by $R\to\infty$ gives the claim.
\end{enumerate}

Recall that if a measurable function on $\R$ admits a Lipschitz
representative (that is, agrees almost everywhere (a.e.) with a
Lipschitz function on $\R$), then this representative is unique and
differentiable a.e. For a function $q:\R\to\R$ that is bounded above,
let $\widehat q$ denote its \emph{least concave majorant}, namely the
pointwise infimum of all concave functions $w\ge q$ on $\R$, and define
the \emph{contact set} of $q$ by
\[
 \mathcal K:=\{z\in\R:\ \widehat q(z)=q(z)\}.
\]

\begin{lemma}\label{lem:contact}
The projection $g_\theta$ has a Lipschitz representative satisfying
$0\le g_\theta\le M/\theta$, and its derivative is given a.e.\ by
\begin{equation}
 g_\theta'(z)=
 \begin{cases}
 1,&m_{\triangle}(z)>\theta,\\
 m(z)/\theta,&m_{\triangle}(z)\le\theta .
 \end{cases}                                                  \label{eq:slope}
\end{equation}
\end{lemma}

\begin{proof}
Note that $v'\le1\iff(v(z)-z)'\le0$, i.e., a primitive of $v(z)-z$ is
concave. Set
\[
 q(z):=\theta^{-1}\int_0^z F(s)\dd s-\frac12 z^2,
\]
so that $q'(z)=\theta^{-1}F(z)-z$ and $q''(z)=\theta^{-1}m(z)-1$.
Since $q(z)=-z^2/2+O(|z|)$ as $|z|\to\infty$, $q$ is bounded above, so
its least concave majorant $\widehat q$ is finite and concave on $\R$.

Set $D:=\widehat q-q\ge0$, so that $\mathcal K=\{D=0\}$. At any
$z\in\mathcal K$, the relation $\widehat q\ge q$ with
$\widehat q(z)=q(z)$ gives
$\widehat q'_-(z)\le q'(z)\le\widehat q'_+(z)$, while concavity of
$\widehat q$ gives $\widehat q'_+(z)\le\widehat q'_-(z)$; hence
$\widehat q'(z)=q'(z)$. It follows that $z\in\mathcal K$ if and only if
the tangent line $w\mapsto q(z)+q'(z)(w-z)$ lies above $q$ on $\R$
(concavity gives ``$\Rightarrow$'', and the minimality of $\widehat q$
gives ``$\Leftarrow$''). By Taylor's formula, for all
$\sigma\in\{-1,1\}$ and $r>0$,
\[
 \begin{aligned}
 q(z+\sigma r)-q(z)-\sigma r q'(z)
   &=\int_0^r(r-s)\bigl(\theta^{-1}m(z+\sigma s)-1\bigr)\dd s\\
   &=\frac{r^2}{2}\bigl(\theta^{-1}U_{\sigma,r}(z)-1\bigr).
 \end{aligned}
\]
Hence
\[
 z\in\mathcal K
 \iff U_{\sigma,r}(z)\le\theta\quad(\forall\,\sigma\in\{-1,1\},\ r>0)
 \iff m_{\triangle}(z)\le\theta,
\]
so $\mathcal K=\{m_{\triangle}\le\theta\}$. Since $m_{\triangle}(z)\to0$
as $|z|\to\infty$, the open set $\Omega:=\{D>0\}=\{m_{\triangle}>\theta\}$
is bounded and $\supp D$ is compact.

Next, we claim that $\widehat q$ is affine on each connected component of
$\Omega$. Fix $z\in\Omega$. By continuity of $\widehat q$ and $D$, we
can choose a closed interval $I_0$ containing $z$ in its interior such
that $\sup_{I_0}\widehat q-\inf_{I_0}\widehat q<\inf_{I_0}D$. Let $\ell$
be the affine function agreeing with $\widehat q$ on $\partial I_0$. By
concavity,
\[
 \inf_{I_0}\widehat q\le\ell\le\widehat q\le\sup_{I_0}\widehat q
 \quad\text{on }I_0,
 \qquad
 \ell\ge\widehat q\ge q
 \quad\text{on }\R\setminus I_0.
\]
Thus $0\le\widehat q-\ell<D$ on $I_0$, so $\ell>q$ on $I_0$, and
$\min\{\widehat q,\ell\}$ is a concave function lying above $q$. By the
minimality of $\widehat q$, we have $\widehat q=\ell$ on $I_0$, which proves the claim.

Since $\Omega$ is open and bounded, each connected component of $\Omega$
is a bounded open interval $(a,b)$ with $a,b\in\mathcal K$. As shown
above, $\widehat q$ is affine on $(a,b)$, and hence on $[a,b]$ by
continuity. Since $\widehat q$ is also
differentiable at $a$ and $b$, its slope on $[a,b]$ equals
$\widehat q'(a)=q'(a)$ and $\widehat q'(b)=q'(b)$. Hence $\widehat q$ is differentiable on $\R$; moreover, $\widehat q'$ is constant on $[a,b]$ for each connected component $(a,b)$
of $\Omega$, and $\widehat q'=q'$ on $\mathcal K$. For any $z_1<z_2$,
if $[z_1,z_2]\cap\mathcal K=\emptyset$, then $[z_1,z_2]$ lies in a
single connected component of $\Omega$, so
$\widehat q'(z_1)=\widehat q'(z_2)$. Otherwise,
$[z_1,z_2]\cap\mathcal K$ is a nonempty compact set, and we define
\[
 \tilde z_1:=\min([z_1,z_2]\cap\mathcal K),\qquad
 \tilde z_2:=\max([z_1,z_2]\cap\mathcal K).
\]
Then $\tilde z_1,\tilde z_2\in\mathcal K$ with
$z_1\le\tilde z_1\le\tilde z_2\le z_2$. Moreover, since
$[z_1,\tilde z_1)$ and $(\tilde z_2,z_2]$ are contained in $\Omega$, we
have $\widehat q'(\tilde z_i)=\widehat q'(z_i)$ for $i=1,2$. Thus
\[
 |\widehat q'(z_2)-\widehat q'(z_1)|
   =|q'(\tilde z_2)-q'(\tilde z_1)|
   \le\|q''\|_{L^\infty_z}(\tilde z_2-\tilde z_1)
   \le\|q''\|_{L^\infty_z}|z_2-z_1|,
\]
so $\widehat q'$ is Lipschitz, $\widehat q\in C^1(\R)$, and
$D\in C_c^1(\R)$.

Define the Lipschitz function $g(z):=z+\widehat q'(z)$. It satisfies
the following properties:
\begin{enumerate}
 \item[(i)] $g'=1$ on $\Omega=\{m_{\triangle}>\theta\}$: this holds
 because $\widehat q'$ is constant on each connected component of
 $\Omega$.
 \item[(ii)] $g'=m/\theta$ a.e.\ on $\mathcal K=\{m_{\triangle}\le\theta\}$:
 since $\widehat q'=q'=F/\theta-z$ on $\mathcal K$, we have $g=F/\theta$
 there, and hence $g'=(F/\theta)'=m/\theta$ a.e.\ on $\mathcal K$.
 \item[(iii)] $g\in\mathcal C_\theta$ and $0\le g\le M/\theta$: by (i),
 (ii), and $0\le m\le m_{\triangle}$, we have $0\le g'\le1$ a.e.; since
 $g=F/\theta$ outside the bounded set $\Omega$, it follows that
 $g\in\mathcal C_\theta$, and monotonicity together with the limits of
 $F/\theta$ at $\pm\infty$ gives $0\le g\le M/\theta$.
\end{enumerate}

Finally, we show that $g=g_\theta$ in $s_\theta+L^2(\R)$. 
Fix $v\in\mathcal C_\theta$. Since $D\in C_c^1(\R)$ is
nonnegative, standard mollification yields nonnegative functions
$\varphi_n\in C_c^\infty(\R)$ with a common compact support such that
$\varphi_n\to D$ and $\varphi_n'\to D'$ uniformly. By the definition of
$\mathcal C_\theta$, $-\int_{\R}v\varphi_n'\dd z\le\int_{\R}\varphi_n\dd z$;
since $v\in L^1_{\mathrm{loc}}(\R)$, letting $n\to\infty$ gives
\[
 -\int_{\R}vD'\dd z\le\int_{\R}D\dd z.
\]

On the other hand, since $g$ is Lipschitz and $D\in C_c^1(\R)$,
integration by parts together with $D=0$ on $\mathcal K$ and $g'=1$ on
$\Omega$ gives
\[
 \int_{\R}gD'\dd z=-\int_{\R}g'D\dd z=-\int_{\R}D\dd z.
\]
Using $F/\theta-g=-D'$, we combine these two estimates to obtain
\[
 \begin{aligned}
 \ip{F/\theta-g}{v-g}_{L^2_z}
   &=-\int_{\R}D'(v-g)\dd z
    =-\int_{\R}vD'\dd z+\int_{\R}gD'\dd z\\
   &\le\int_{\R}D\dd z-\int_{\R}D\dd z=0.
 \end{aligned}
\]
Since $g\in\mathcal C_\theta$ by (iii) and $v\in\mathcal C_\theta$ is
arbitrary, the variational characterization \eqref{eq:var} implies
$g=g_\theta$ in $s_\theta+L^2(\R)$. Thus $g$ is the Lipschitz
representative of $g_\theta$, and the lemma follows from (i)--(iii).
\end{proof}

By Lemma~\ref{lem:contact}, the function
$m_{\triangle}$ determines the contact set for every level $\theta>0$. In
Subsection~\ref{subsec:log}, we will average the slope \eqref{eq:slope}
logarithmically in $\theta$ (see \eqref{eq:hexplicit}), which produces
a $\log m_{\triangle}$ term and hence the ratio
$m_{\triangle}'/m_{\triangle}$ upon differentiation in $z$. Applying
Cauchy--Schwarz to the resulting error term leads to the integral
$\int_{\R}(|m_{\triangle}'|^2/m_{\triangle})\dd z$. The next lemma
establishes a Fisher-type estimate bounding this integral in terms of
$\|(\sqrt m)'\|_{L^2_z}^2$.

\begin{lemma}\label{lem:fisher}
There exists an absolute constant $C_0>0$ such that
\begin{equation}
 \int_{\R}\frac{|m_{\triangle}'|^2}{m_{\triangle}}\dd z
     \le C_0\int_{\R}|(\sqrt m)'|^2\dd z.                    \label{eq:fisher}
\end{equation}
\end{lemma}

\begin{proof}
For every $\sigma\in\{-1,1\}$, $r>0$, and $z_1,z_2\in\R$,
\eqref{eq:triangular} gives
\[
 |U_{\sigma,r}(z_1)-U_{\sigma,r}(z_2)|
   \le\frac{2}{r^2}\int_0^r(r-s)\|m'\|_{L^\infty_z}|z_1-z_2|\dd s
   =\|m'\|_{L^\infty_z}|z_1-z_2|.
\]
Thus $m_{\triangle}=\sup_{\sigma,r}U_{\sigma,r}$ is Lipschitz with
constant at most $\|m'\|_{L^\infty_z}$, and in particular is
differentiable a.e.\ on $\R$.

Set $f:=\sqrt m$. Since $m$ is continuous and $f\in H^1(\R)$, $f$ is
locally absolutely continuous on $\R$. Fix a point $z\in\R$ where
$m_{\triangle}$ is differentiable. If $m_{\triangle}(z)=m(z)$, then
$m(z)=m_{\triangle}(z)>0$, and the nonnegative function
$m_{\triangle}-m$ vanishes at $z$, so it attains a minimum there. Hence
$m_{\triangle}'(z)=m'(z)$, which gives
\[
 \frac{|m_{\triangle}'(z)|}{\sqrt{m_{\triangle}(z)}}
   =\frac{|m'(z)|}{\sqrt{m(z)}}=2|f'(z)|.
\]

Now suppose $m_{\triangle}(z)>m(z)$. Since $U_{\sigma,r}(z)$ is
continuous in $r\in(0,\infty)$ with
$\lim_{r\downarrow0}U_{\sigma,r}(z)=m(z)$ and
$\lim_{r\to\infty}U_{\sigma,r}(z)=0$, the supremum in
\eqref{eq:triangular} is attained at some pair $(\sigma_0,r_0)$ with
$\sigma_0\in\{-1,1\}$ and $0<r_0<\infty$. For any $\sigma\in\{-1,1\}$
and $r>0$, denote the average of $m$ over the closed interval with
endpoints $z$ and $z+\sigma r$ by
\[
 \overline m_{\sigma,r}(z):=\frac1r\int_0^r m(z+\sigma s)\dd s.
\]
Differentiating \eqref{eq:triangular} in $r$, and differentiating in
$z$ followed by integration by parts in $s$, we obtain
\[
 \partial_r U_{\sigma,r}(z)
  =\frac2r\bigl(\overline m_{\sigma,r}(z)-U_{\sigma,r}(z)\bigr),
 \qquad
 \partial_z U_{\sigma,r}(z)
  =\frac{2\sigma}{r}\bigl(\overline m_{\sigma,r}(z)-m(z)\bigr).
\]
At $(\sigma,r)=(\sigma_0,r_0)$, we have $\partial_r U_{\sigma,r}(z)=0$,
so $\overline m_{\sigma_0,r_0}(z)=U_{\sigma_0,r_0}(z)=m_{\triangle}(z)$.
Moreover, the nonnegative function
$m_{\triangle}-U_{\sigma_0,r_0}$ is differentiable at $z$
and attains a minimum there, so
$m_{\triangle}'(z)=\partial_z U_{\sigma_0,r_0}(z)$. Evaluating the
second identity at $(\sigma_0,r_0)$ gives
\[
 m_{\triangle}'(z)=\frac{2\sigma_0}{r_0}\bigl(m_{\triangle}(z)-m(z)\bigr).
\]

Let $I_{r_0}$ be the closed interval with endpoints $z$ and
$z+\sigma_0 r_0$. Since
$m_{\triangle}(z)=U_{\sigma_0,r_0}(z)\le\sup_{I_{r_0}}m$, we have
$f(z)\le\sqrt{m_{\triangle}(z)}\le\sup_{I_{r_0}}f$. Using $m=f^2$, we obtain
\[
 \begin{aligned}
 \frac{|m_{\triangle}'(z)|}{\sqrt{m_{\triangle}(z)}}
   &=\frac{2}{r_0}\frac{\sqrt{m_{\triangle}(z)}+f(z)}{\sqrt{m_{\triangle}(z)}}
     \bigl(\sqrt{m_{\triangle}(z)}-f(z)\bigr)\\
   &\le\frac{4}{r_0}\bigl(\sqrt{m_{\triangle}(z)}-f(z)\bigr)
    \le\frac{4}{r_0}\bigl(\sup_{I_{r_0}}f-f(z)\bigr)\\
   &\le\frac{4}{r_0}\int_{I_{r_0}} |f'(s)|\dd s
    \le4\HL(|f'|)(z),
 \end{aligned}
\]
where $\HL$ is the uncentered Hardy--Littlewood maximal operator on
$\R$, defined for $w\in L^1_{\mathrm{loc}}(\R)$ by
\[
 \HL w(z):=\sup_{B\ni z}\frac{1}{|B|}\int_B|w(s)|\dd s,
\]
with the supremum taken over all bounded intervals $B$ of positive
length containing $z$.

Combining the two cases, we have
\[
 \frac{|m_{\triangle}'(z)|}{\sqrt{m_{\triangle}(z)}}
   \le 2|f'(z)|+4\HL(|f'|)(z)
 \qquad\text{for a.e.\ }z\in\R.
\]
Squaring both sides, integrating over $\R$, and using the $L^2_z$
boundedness of the Hardy--Littlewood maximal operator,
$\|\HL w\|_{L^2_z}\le C_{\mathrm{HL}}\|w\|_{L^2_z}$, we obtain \eqref{eq:fisher}.
\end{proof}

\subsection{Logarithmic averaging}\label{subsec:log}
For a fixed level $\theta>0$, the projection $g_\theta$ from
Subsection~\ref{subsec:proj} is not yet suitable as a virial multiplier:
at a boundary point of $\{m_{\triangle}>\theta\}$ we may have
$m<\theta=m_{\triangle}$, so the two branches of $g_\theta'$ in
\eqref{eq:slope} can meet with a jump. In the directional virial
identity of Section~\ref{sec:virial} (see \eqref{eq:jk}), however, the
derivative of the slope is paired with $m'$, which requires a Lipschitz
slope. To remove this jump, we average $g_\theta$ logarithmically in
$\theta$.

Let $h_\theta(z)$ denote the right-hand side of \eqref{eq:slope}. Since
$m$ and $m_{\triangle}$ are continuous, $h_\theta(z)$ is jointly
measurable in $(\theta,z)$ and takes values in $[0,1]$. By
Lemma~\ref{lem:contact}, for each $\theta>0$ the Lipschitz
representative $g_\theta$ satisfies $g_\theta'=h_\theta$ a.e.\ and
$\lim_{z\to-\infty}g_\theta(z)=0$, so
$g_\theta(z)=\int_{-\infty}^z h_\theta(s)\dd s$ is also jointly
measurable in $(\theta,z)$. Fix $0<\theta_0<\theta_1$ with
$L:=\log(\theta_1/\theta_0)>1$, let
$\dd\mu(\theta):=\dd\theta/(L\theta)$ be the uniform probability measure
in $\log\theta$ on $[\theta_0,\theta_1]$, and define
\begin{equation}
 G(z):=\int_{\theta_0}^{\theta_1}g_\theta(z)\dd\mu(\theta)
      =\frac1L\int_{\theta_0}^{\theta_1}g_\theta(z)\frac{\dd\theta}{\theta}.
                                                                \label{eq:G}
\end{equation}
The next lemma shows that averaging smooths the slope to a Lipschitz
function and yields an $O(L^{-1})$ bound on the derivative error term.

\begin{lemma}\label{lem:log}
The function $G$ defined in \eqref{eq:G} is Lipschitz with
$0\le G\le M/\theta_0$. Its derivative $h:=G'$ has a Lipschitz
representative satisfying $0\le h\le1$, $h=1$ on $\{m\ge\theta_1\}$, and
\begin{equation}
 \left|\int_{\R} h'm'\dd z\right|
       \le\frac{C_1}{L}\int_{\R}|(\sqrt m)'|^2\dd z             \label{eq:error}
\end{equation}
for an absolute constant $C_1>0$.
\end{lemma}

\begin{proof}
By Lemma~\ref{lem:contact}, we have
$0\le g_\theta\le M/\theta\le M/\theta_0$ and $0\le h_\theta\le1$ for
all $\theta\in[\theta_0,\theta_1]$. Averaging with respect to the
probability measure $\mu$ immediately gives $0\le G\le M/\theta_0$.
Moreover, since $g_\theta(z)=\int_{-\infty}^z h_\theta(s)\dd s$,
Tonelli's theorem gives, for any $a<b$,
\[
 G(b)-G(a)=\int_a^b\int_{\theta_0}^{\theta_1}
               h_\theta(s)\dd\mu(\theta)\dd s .
\]
Hence $G$ is Lipschitz on $\R$, and its derivative $h=G'$ is given for
a.e.\ $z\in\R$ by
\[
 h(z)=\int_{\theta_0}^{\theta_1}h_\theta(z)\dd\mu(\theta)
     =\frac1L\int_{\theta_0}^{\theta_1}h_\theta(z)\frac{\dd\theta}{\theta}.
\]
If $m_{\triangle}(z)\le\theta_0$, then $h_\theta(z)=m(z)/\theta$ on
$[\theta_0,\theta_1]$. If $m_{\triangle}(z)\ge\theta_1$, then $h_\theta(z)=1$ on
$[\theta_0,\theta_1)$. Finally, for
$\theta_0<m_{\triangle}(z)<\theta_1$, splitting the integral at
$m_{\triangle}(z)$ gives
\[
 \begin{aligned}
  Lh(z)&=\int_{\theta_0}^{m_{\triangle}(z)}\frac{\dd\theta}{\theta}
        +m(z)\int_{m_{\triangle}(z)}^{\theta_1}\frac{\dd\theta}{\theta^2}\\
       &=\log\bigl(m_{\triangle}(z)/\theta_0\bigr)
        +m(z)\bigl(m_{\triangle}(z)^{-1}-\theta_1^{-1}\bigr).
 \end{aligned}
\]
Combining the three cases yields the continuous representative of $h$
given by
\begin{equation}
 h(z)=
 \begin{cases}
 \frac{m}{L}(\theta_0^{-1}-\theta_1^{-1}), &m_{\triangle}\le\theta_0,\\
 \frac1L\bigl(\log(m_{\triangle}/\theta_0)+m(m_{\triangle}^{-1}-\theta_1^{-1})\bigr), &\theta_0<m_{\triangle}<\theta_1,\\
 1, &m_{\triangle}\ge\theta_1 .
 \end{cases}                                                 \label{eq:hexplicit}
\end{equation}
Since $0\le h_\theta\le1$ and $\mu$ is a probability measure, this
representative satisfies $0\le h\le1$ on $\R$; since
$m_{\triangle}\ge m$, the third branch of \eqref{eq:hexplicit} also
shows that $h=1$ on $\{m\ge\theta_1\}$.

Note that the adjacent branches in \eqref{eq:hexplicit} agree on
$\{m_{\triangle}=\theta_0\}$ and $\{m_{\triangle}=\theta_1\}$. To
differentiate across these level sets, we set
\[
 \tau:=\min\{\theta_1,\max\{m_{\triangle},\theta_0\}\}.
\]
Then $\tau$ is Lipschitz with $\theta_0\le\tau\le\theta_1$, and
\eqref{eq:hexplicit} can be rewritten as
$Lh=\log(\tau/\theta_0)+m(\tau^{-1}-\theta_1^{-1})$; in particular,
$h$ is Lipschitz on $\R$. Since
$\{\tau=\theta_0\}=\{m_{\triangle}\le\theta_0\}$ and
$\{\tau=\theta_1\}=\{m_{\triangle}\ge\theta_1\}$, and a Lipschitz
function has derivative zero a.e.\ on each of its level sets, we have
$\tau'=0$ a.e.\ on
$\{m_{\triangle}\le\theta_0\}\cup\{m_{\triangle}\ge\theta_1\}$. On the
open set $\{\theta_0<m_{\triangle}<\theta_1\}$ we have
$\tau=m_{\triangle}$ and hence $\tau'=m_{\triangle}'$ a.e. By the chain
rule,
\[
 Lh'=\left(1-\frac{m}{\tau}\right)\frac{\tau'}{\tau}
    +(\tau^{-1}-\theta_1^{-1})m'
 \qquad\text{a.e.\ on }\R.
\]
On $\{\theta_0<m_{\triangle}<\theta_1\}$, we have
$m\le m_{\triangle}=\tau$, so $0\le 1-m/\tau\le1$ and
$0\le\tau^{-1}-\theta_1^{-1}\le m_{\triangle}^{-1}$. On
$\{m_{\triangle}\le\theta_0\}$, we have $\tau'=0$ a.e.\ and
$0\le\tau^{-1}-\theta_1^{-1}\le\theta_0^{-1}\le m_{\triangle}^{-1}$. On
$\{m_{\triangle}\ge\theta_1\}$, we have $\tau=\theta_1$ and $\tau'=0$
a.e., so $h'=0$ a.e. Combining all three regions gives
\[
 L|h'|\le\frac{|m'|+|m_{\triangle}'|}{m_{\triangle}}
 \qquad\text{a.e.\ on }\R.
\]

Multiplying by $|m'|$, integrating over $\R$, and applying the
Cauchy--Schwarz inequality, we obtain
\[
 L\int_{\R}|h'm'|\dd z
 \le \int_{\R}\frac{|m'|^2}{m_{\triangle}}\dd z
   +\left(\int_{\R}\frac{|m'|^2}{m_{\triangle}}\dd z\right)^{1/2}
    \left(\int_{\R}\frac{|m_{\triangle}'|^2}{m_{\triangle}}\dd z\right)^{1/2}.
\]
Set $f:=\sqrt m$. Since $m'=2ff'$ a.e.\ and $m\le m_{\triangle}$, we
have $|m'|^2/m_{\triangle}=4(m/m_{\triangle})|f'|^2\le 4|f'|^2$ a.e.
Using Lemma~\ref{lem:fisher} for the last term, we obtain
\eqref{eq:error} with $C_1=4+2\sqrt{C_0}$.
\end{proof}

\section{Virial inequalities for Schwartz solutions}\label{sec:virial}

In this section, we apply the averaged multiplier $G$ from
Section~\ref{sec:multiplier} to a nonzero Schwartz solution
$u\in C(I;\Sch(\R^3))$ of \eqref{eq:nls}. Fix a coordinate direction
$k\in\{1,2,3\}$, write $z=x_k$, denote the
transverse coordinates by $y\in\R^2$, and identify $x\in\R^3$ with
$(z,y)\in\R\times\R^2$. Define the planar mass and momentum densities
\begin{align}
 m(z,t)&=\int_{\R^2}|u(z,y,t)|^2\dd y, &
 p(z,t)&=\im\int_{\R^2}\bar u(z,y,t)\partial_z u(z,y,t)\dd y .
                                                               \label{eq:mp}
\end{align}
Since $\sqrt{m(z,t)}=\|u(z,\cdot,t)\|_{L^2_y(\R^2)}$, for any $\eta>0$,
by the reverse triangle inequality in $L^2_y(\R^2)$, the fundamental
theorem of calculus, and Cauchy--Schwarz, we have
\[
 \begin{aligned}
  &\int_{\R}\left|\frac{\sqrt{m(z+\eta,t)}-\sqrt{m(z,t)}}{\eta}\right|^2\dd z\\
  &\qquad\le\int_{\R}\!\int_{\R^2}
        \left|\frac{u(z+\eta,y,t)-u(z,y,t)}{\eta}\right|^2\dd y\dd z\\
  &\qquad\le\frac{1}{\eta}\int_0^\eta\int_{\R}\!\int_{\R^2}
        |\partial_z u(z+s,y,t)|^2\dd y\dd z\dd s
       =\|\partial_z u(t)\|_{L^2_x}^2.
 \end{aligned}
\]
Combining this estimate with $\sqrt{m(\cdot,t)}\in L^2(\R)$, we obtain
$\sqrt{m(\cdot,t)}\in H^1(\R)$ and
\begin{equation}
 \int_{\R}|\partial_z\sqrt{m(z,t)}|^2\dd z
       \le\|\partial_z u(t)\|_{L^2_x}^2.                            \label{eq:marggrad}
\end{equation}
For each $t\in I$, $m(\cdot,t)$ is a nonnegative Schwartz function with
conserved total mass $M(u)>0$ and $\sqrt{m(\cdot,t)}\in H^1(\R)$. Hence the construction of Section~\ref{sec:multiplier} applies at each
$t\in I$, yielding $g_\theta(z,t)$ and the averaged multiplier $G(z,t)$
with slope $h=\partial_z G$. In particular, combining \eqref{eq:error}
with \eqref{eq:marggrad} gives
\begin{equation}
 \left|\int_{\R}\partial_z h\,\partial_z m\dd z\right|
   \le\frac{C_1}{L}\|\partial_z u(t)\|_{L^2_x}^2,              \label{eq:spatial_err}
\end{equation}
which will control the derivative error term in the directional
virial identity \eqref{eq:jk}.

\subsection{Time dependence along the NLS flow}\label{subsec:time}
We first study the time dependence of $G$ in order to justify
differentiating $\int_{\R}Gp\dd z$ in time and show that
$\int_{\R}G_t p\dd z\le0$. To this end, we relate $p$ to the cumulative
mass
$F(z,t)=\int_{-\infty}^z m(s,t)\dd s$. Multiplying \eqref{eq:nls} by
$2\bar u$ and taking imaginary parts gives the continuity equation
\[
 \partial_t(|u|^2)=-2\nabla\cdot\im(\bar u\nabla u)
                  =-2\partial_z\im(\bar u\,\partial_z u)
                   -2\nabla_y\cdot\im(\bar u\nabla_y u).
\]
Since $u\in C(I;\Sch(\R^3))$ is a Schwartz solution, we may
differentiate under the integral sign. Integrating the above identity
over $y\in\R^2$ and noting that the transverse divergence term vanishes,
we obtain $m_t(z,t)=-2\partial_z p(z,t)$. Differentiating $F(z,t)$ in $t$ and
using $\lim_{s\to-\infty}p(s,t)=0$, we obtain
\begin{equation}
 F_t(z,t)=\int_{-\infty}^z m_t(s,t)\dd s
         =-2\int_{-\infty}^z \partial_s p(s,t)\dd s
         =-2p(z,t).                                           \label{eq:continuity}
\end{equation}
Moreover, by mass conservation, $\lim_{z\to+\infty}F(z,t)=M(u)$ is
independent of $t$, so for each $\theta>0$ the convex set
$\mathcal C_\theta$ does not depend on $t$. Comparing the projections
$g_\theta(\cdot,t)$ at different times via the variational inequality
\eqref{eq:var} therefore yields the following estimate.

\begin{lemma}\label{lem:time}
The averaged multiplier $G(z,t)$ constructed from \eqref{eq:mp} is a locally
Lipschitz curve in $G(0)+L^2(\R)$, where $G(t):=G(\cdot,t)$, and for
a.e.\ $t\in I$,
\begin{equation}
 \int_{\R}G_t(z,t)p(z,t)\dd z\le0 .                          \label{eq:good}
\end{equation}
\end{lemma}

\begin{proof}
Write $F(t):=F(\cdot,t)$, and use the same notation for $p(t)$ and
$g_\theta(t)$. Since $u\in C(I;\Sch(\R^3))$, the map
$t\mapsto p(t)$ is continuous from $I$ into $L^2(\R)$. Integrating
\eqref{eq:continuity} in time gives
$F(t)-F(s)=-2\int_s^t p(r)\dd r$ in $L^2(\R)$, so $t\mapsto F(t)-F(0)$
is a $C^1$ curve in $L^2(\R)$ with derivative $F_t=-2p$.

Fix $s,t\in I$. Since both $g_\theta(s)$
and $g_\theta(t)$ belong to the convex set $\mathcal C_\theta$,
applying \eqref{eq:var} at time $t$ with $v=g_\theta(s)$ and at time $s$
with $v=g_\theta(t)$ gives
\[
 \ip{F(t)/\theta-g_\theta(t)}{g_\theta(s)-g_\theta(t)}_{L^2_z}\le0,
 \qquad
 \ip{F(s)/\theta-g_\theta(s)}{g_\theta(t)-g_\theta(s)}_{L^2_z}\le0.
\]
Adding these two inequalities and multiplying by $\theta$, we obtain
\begin{equation}
 \ip{g_\theta(t)-g_\theta(s)}{F(t)-F(s)}_{L^2_z}
      \ge\theta\|g_\theta(t)-g_\theta(s)\|_{L^2_z}^2 .       \label{eq:mono_theta}
\end{equation}
By Cauchy--Schwarz, it follows that
$\|g_\theta(t)-g_\theta(s)\|_{L^2_z}\le\theta^{-1}\|F(t)-F(s)\|_{L^2_z}$.

For fixed $s,t\in I$, the difference $g_\theta(z,t)-g_\theta(z,s)$ is
jointly measurable in $(\theta,z)$. By Cauchy--Schwarz in $z$ and the
above estimate,
\[
 \begin{aligned}
 &\int_{\theta_0}^{\theta_1}\!\int_{\R}
   \bigl|\bigl(g_\theta(z,t)-g_\theta(z,s)\bigr)
         \bigl(F(z,t)-F(z,s)\bigr)\bigr|\dd z\dd\mu(\theta)\\
 &\qquad\le\int_{\theta_0}^{\theta_1}
   \|g_\theta(t)-g_\theta(s)\|_{L^2_z}\|F(t)-F(s)\|_{L^2_z}\dd\mu(\theta)\\
 &\qquad\le\theta_0^{-1}\|F(t)-F(s)\|_{L^2_z}^2<\infty .
 \end{aligned}
\]
Hence Fubini's theorem and \eqref{eq:G} yield
\[
 \int_{\theta_0}^{\theta_1}
   \ip{g_\theta(t)-g_\theta(s)}{F(t)-F(s)}_{L^2_z}\dd\mu(\theta)
 =\ip{G(t)-G(s)}{F(t)-F(s)}_{L^2_z}.
\]
On the other hand, since $\mu$ is a probability measure, Cauchy--Schwarz
in $\theta$ gives
\[
 |G(z,t)-G(z,s)|^2
   \le\int_{\theta_0}^{\theta_1}|g_\theta(z,t)-g_\theta(z,s)|^2\dd\mu(\theta)
\]
for every $z\in\R$. Integrating in $z$ and applying Tonelli's theorem,
we see that $G(t)-G(s)\in L^2(\R)$ with
\[
 \|G(t)-G(s)\|_{L^2_z}^2
 \le\int_{\theta_0}^{\theta_1}\|g_\theta(t)-g_\theta(s)\|_{L^2_z}^2\dd\mu(\theta).
\]
Averaging \eqref{eq:mono_theta} with respect to $\mu$ and using
$\theta\ge\theta_0$, we obtain
\begin{equation}
 \begin{aligned}
  \ip{G(t)-G(s)}{F(t)-F(s)}_{L^2_z}
  &\ge\int_{\theta_0}^{\theta_1}\theta\|g_\theta(t)-g_\theta(s)\|_{L^2_z}^2\dd\mu(\theta)\\
  &\ge\theta_0\|G(t)-G(s)\|_{L^2_z}^2.
 \end{aligned}                                                \label{eq:mono_G}
\end{equation}

In particular, Cauchy--Schwarz gives
\[
 \|G(t)-G(s)\|_{L^2_z}\le\theta_0^{-1}\|F(t)-F(s)\|_{L^2_z}.
\]
Since
$t\mapsto F(t)-F(0)$ is $C^1$ in $L^2(\R)$, it follows that
$t\mapsto G(t)-G(0)$ is locally Lipschitz with values in the Hilbert
space $L^2(\R)$, and hence strongly differentiable in $L^2(\R)$ for
a.e.\ $t\in I$. At any such time $t$, dividing \eqref{eq:mono_G} by
$(t-s)^2$ and letting $s\to t$, both difference quotients converge in
$L^2(\R)$, giving $\ip{G_t}{F_t}_{L^2_z}\ge\theta_0\|G_t\|_{L^2_z}^2$.
Combining this with \eqref{eq:continuity}, we obtain
\[
 \int_{\R}G_t(z,t)p(z,t)\dd z
   =-\frac12\ip{G_t}{F_t}_{L^2_z}
   \le-\frac{\theta_0}{2}\|G_t\|_{L^2_z}^2\le0,
\]
which proves \eqref{eq:good}.
\end{proof}

\subsection{The directional virial identity}\label{subsec:jk}
Define the planar integrals
\[
 a(z,t)=\int_{\R^2}|\partial_z u(z,y,t)|^2\dd y,\qquad
 b(z,t)=\int_{\R^2}|u(z,y,t)|^4\dd y.
\]
Using \eqref{eq:nls} in the form $\partial_t u=i(\Delta u+|u|^2 u)$, we compute
\[
 \partial_t\im(\bar u\,\partial_z u)
  =\re(\bar u\,\partial_z\Delta u-\Delta\bar u\,\partial_z u)
   +\frac12\partial_z(|u|^4).
\]
For the first term on the right-hand side, the product rule gives
\[
 \re(\bar u\,\partial_z\Delta u-\Delta\bar u\,\partial_z u)
  =\partial_z\re(\bar u\Delta u)-2\re(\Delta\bar u\,\partial_z u),
\]
while
$\nabla\cdot\re(\nabla\bar u\,\partial_z u)=\re(\Delta\bar u\,\partial_z u)+\frac12\partial_z(|\nabla u|^2)$
and $\frac12\Delta(|u|^2)=\re(\bar u\Delta u)+|\nabla u|^2$. Combining
these identities yields the local momentum law
\[
 \partial_t\im(\bar u\,\partial_z u)
  =-2\nabla\cdot\re(\nabla\bar u\,\partial_z u)
   +\frac12\partial_z\Delta(|u|^2)
   +\frac12\partial_z(|u|^4).
\]
Splitting $\nabla=(\partial_z,\nabla_y)$ and
$\Delta=\partial_z^2+\Delta_y$, we integrate over $y\in\R^2$. Since
$u\in C(I;\Sch(\R^3))$ is a Schwartz solution, we may differentiate
under the integral sign and the transverse divergence terms vanish,
giving
\begin{equation}
 p_t=-2\partial_z a+\frac12\partial_z^3 m+\frac12\partial_z b. \label{eq:momentum}
\end{equation}

For this fixed direction $k$, define the directional virial quantity
$J_k(t)=\int_{\R}G(z,t)p(z,t)\dd z$. For any fixed $t_0\in I$, we split
\[
 J_k(t)=\int_{\R}G(z,t_0)p(z,t)\dd z+\ip{G(t)-G(t_0)}{p(t)}_{L^2_z}.
\]
Here $G(t_0)$ is bounded, while $t\mapsto G(t)-G(t_0)$ is locally
Lipschitz with values in $L^2(\R)$ by Lemma~\ref{lem:time}. Since $u$
is a Schwartz solution, $t\mapsto p(t)$ is $C^1$ with values in
$L^1(\R)\cap L^2(\R)$. Hence $J_k$ is locally absolutely continuous
on $I$, and the product rule holds for a.e.\ $t\in I$, giving
\[
 J_k'(t)=\int_{\R}G_t p\dd z+\int_{\R}G p_t\dd z.
\]
Substituting \eqref{eq:momentum} into the second integral and
integrating by parts in $z$, using
\[
 \int_{\R}G\,\partial_z^3 m\dd z
   =-\int_{\R}h\,\partial_z^2 m\dd z
   =\int_{\R}\partial_z h\,\partial_z m\dd z,
\]
we obtain
\begin{equation}
 J_k'(t)=\int_{\R}G_t p\dd z+2\int_{\R}ha\dd z-\frac12\int_{\R}hb\dd z
              +\frac12\int_{\R}\partial_z h\,\partial_z m\dd z. \label{eq:jk}
\end{equation}
These integrations by parts are justified in the weak sense because $G$
and $h$ are bounded and Lipschitz by Lemma~\ref{lem:log}, while $a$,
$b$, $m$, and their derivatives decay rapidly.

Since $0\le h\le 1$, the second term on the right-hand side of
\eqref{eq:jk} is bounded above by $2\int_{\R} a\dd z$. Comparing the
third term $-\frac12\int_{\R}hb\dd z$ with the full $L^4$ term
$-\frac12\int_{\R}b\dd z=-\frac12\|u(t)\|_{L^4_x}^4$ leaves the error
$\frac12\int_{\R}(1-h)b\dd z$. To estimate it, we apply the
two-dimensional Gagliardo--Nirenberg inequality on each plane $x_k=z$,
obtaining
\[
 b(z,t)\le C\,m(z,t)
       \int_{\R^2}|\nabla_y u(z,y,t)|^2\dd y
\]
for an absolute constant $C>0$.
Since $h=1$ whenever $m\ge\theta_1$ and $0\le h\le1$, it follows that
\begin{equation}
 \int_{\R}(1-h)b\dd z
     \le C\theta_1\|\nabla_y u(t)\|_{L^2_x}^2 . \label{eq:slice}
\end{equation}
Thus $\theta_1$ alone controls the $L^4$ loss. With $\theta_1$
fixed, we can still decrease $\theta_0$ and make
$L=\log(\theta_1/\theta_0)$ arbitrarily large in the derivative error
estimate \eqref{eq:spatial_err}.

\begin{proposition}\label{prop:virial}
Let $u$ be a Schwartz solution of \eqref{eq:nls} on an interval $I$
containing $0$ satisfying \eqref{eq:hyp} (and hence \eqref{eq:trap}),
with $\delta\in(0,1]$ defined in Theorem~\ref{thm:main}. There exist
parameters $\theta_1>0$ and $L>1$, and a constant $c_\delta>0$,
depending only on $\delta$, with $\delta\mapsto c_\delta$ continuous on
$(0,1]$, such that when $\theta_0=\theta_1e^{-L}$, the sum
$J=\sum_{k=1}^3 J_k$ is locally absolutely continuous on $I$ and
satisfies
\begin{equation}
 J'(t)\le-\frac12\delta\|\nabla u(t)\|_{L^2_x}^2,\qquad
 |J(t)|\le c_\delta M(u)^{3/2}\|\nabla u(t)\|_{L^2_x},            \label{eq:goal}
\end{equation}
where the first inequality holds for a.e.\ $t\in I$ and the second holds
for all $t\in I$.
\end{proposition}

\begin{proof}
Since $J_k$ is locally absolutely continuous for each $k\in\{1,2,3\}$,
the sum $J=\sum_{k=1}^3 J_k$ is also locally absolutely continuous. Fix a coordinate direction $k\in\{1,2,3\}$,
with $z=x_k$ and $y\in\R^2$ as above, and write $\partial_j=\partial_{x_j}$. Then
\[
 \int_{\R}a\dd z=\|\partial_z u(t)\|_{L^2_x}^2=\|\partial_k u(t)\|_{L^2_x}^2,
 \qquad
 \int_{\R}b\dd z=\|u(t)\|_{L^4_x}^4,
\]
and $\|\nabla_y u(t)\|_{L^2_x}^2=\sum_{j\ne k}\|\partial_j u(t)\|_{L^2_x}^2$.
In \eqref{eq:jk}, we drop the nonpositive first term by
Lemma~\ref{lem:time} and bound the second using $h\le1$. For the
remaining two terms, we write
$-\int_{\R}hb\dd z=-\int_{\R}b\dd z+\int_{\R}(1-h)b\dd z$ and apply
\eqref{eq:slice} and \eqref{eq:spatial_err}. This yields, for a.e.\
$t\in I$,
\[
 J_k'(t)\le 2\|\partial_k u(t)\|_{L^2_x}^2
           -\frac12\|u(t)\|_{L^4_x}^4
           +\frac12 C\theta_1
            \sum_{j\ne k}\|\partial_j u(t)\|_{L^2_x}^2
           +\frac{C_1}{2L}\|\partial_k u(t)\|_{L^2_x}^2.
\]
Summing over $k=1,2,3$ and using
$\sum_{k=1}^3\sum_{j\ne k}\|\partial_j u(t)\|_{L^2_x}^2=2\|\nabla u(t)\|_{L^2_x}^2$
together with $2K(u)=2\|\nabla u\|_{L^2_x}^2-\frac32\|u\|_{L^4_x}^4$, we
obtain
\[
 J'(t)\le 2K(u(t))+\left(C\theta_1
                +\frac{C_1}{2L}\right)\|\nabla u(t)\|_{L^2_x}^2 .
\]
We set $\theta_1=\delta/(4C)$, $L=\max\{2,2C_1/\delta\}>1$, and
$\theta_0=\theta_1e^{-L}$. Then $C\theta_1=\delta/4$ and
$C_1/(2L)\le\delta/4$. Combining this with \eqref{eq:trap} gives the
first inequality in \eqref{eq:goal}.

For the second inequality in \eqref{eq:goal}, fix $k\in\{1,2,3\}$. Using
\eqref{eq:mp} and Cauchy--Schwarz, we have
\[
 \begin{aligned}
 \int_{\R}|p(z,t)|\dd z
 &\le\int_{\R}\!\int_{\R^2}|u(z,y,t)|\,|\partial_z u(z,y,t)|\dd y\dd z\\
 &=\int_{\R^3}|u(x,t)|\,|\partial_k u(x,t)|\dd x
  \le M(u)^{1/2}\|\partial_k u(t)\|_{L^2_x},
 \end{aligned}
\]
while Lemma~\ref{lem:log} gives $\|G(t)\|_{L^\infty_z}\le M(u)/\theta_0$.
Hence
\[
 |J_k(t)|
 \le\|G(t)\|_{L^\infty_z}\int_{\R}|p(z,t)|\dd z
 \le\frac{1}{\theta_0}M(u)^{3/2}\|\partial_k u(t)\|_{L^2_x}.
\]
Summing over $k=1,2,3$ and using
$\sum_{k=1}^3\|\partial_k u\|_{L^2_x}\le\sqrt3\|\nabla u\|_{L^2_x}$, we
obtain
\[
 |J(t)|
 \le\sum_{k=1}^3|J_k(t)|
 \le\frac{\sqrt3}{\theta_0}M(u)^{3/2}\|\nabla u(t)\|_{L^2_x} ,
\]
which gives the second inequality in \eqref{eq:goal} with
\[
 c_\delta=\frac{\sqrt3}{\theta_0}
         =\frac{4\sqrt3\,C}{\delta}\exp\bigl(\max\{2,2C_1/\delta\}\bigr).
\]
Clearly $c_\delta$ depends only on $\delta$ and is continuous on
$(0,1]$.
\end{proof}

\section{Proof of Theorem~\ref{thm:main}}\label{sec:lifespan}

\begin{proof}
We first assume that $u\in C([t_1,t_2];\Sch(\R^3))$ is a Schwartz
solution of \eqref{eq:nls} satisfying \eqref{eq:hyp} on a compact
interval $[t_1,t_2]$ containing $0$. We prove an upper bound on
$t_2-t_1$. Let $\delta\in(0,1]$ be as in Theorem~\ref{thm:main}; then
$u$ satisfies \eqref{eq:trap} on $[t_1,t_2]$. By \eqref{eq:trap} and
Proposition~\ref{prop:virial}, we have
$\|\nabla u(t)\|_{L^2_x}^2\ge 6M(Q)E(Q)/M(u)$ and
$\|\nabla u(t)\|_{L^2_x}^2\ge J(t)^2/(c_\delta^2 M(u)^3)$ for all
$t\in[t_1,t_2]$. Averaging these two lower bounds and substituting into
the first inequality of \eqref{eq:goal},
we obtain, for a.e.\ $t\in[t_1,t_2]$,
\[
 -J'(t)\ge\frac{\delta}{4}
     \left(\frac{6M(Q)E(Q)}{M(u)}+\frac{J(t)^2}{c_\delta^2 M(u)^3}\right).
\]
By the Pohozaev identity
$6M(Q)E(Q)=\|Q\|_{L^2_x}^2\|\nabla Q\|_{L^2_x}^2$, setting
$Y(t)=J(t)/(c_\delta M(u)\Qprod)$ yields
\[
 -Y'(t)\ge\frac{\delta\Qprod}{4c_\delta M(u)^2}\bigl(1+Y(t)^2\bigr),
\]
and hence
\[
 -\frac{d}{dt}\arctan Y(t)
     =-\frac{Y'(t)}{1+Y(t)^2}
     \ge\frac{\delta\Qprod}{4c_\delta M(u)^2}
      \qquad\text{for a.e.\ }t\in[t_1,t_2].
\]
Since $Y$ is absolutely continuous on $[t_1,t_2]$, the function
$\arctan Y$ is also absolutely continuous on $[t_1,t_2]$.
Integrating over $[t_1,t_2]$ and using $|\arctan Y|<\pi/2$ gives
$\arctan Y(t_1)-\arctan Y(t_2)\le\pi$, and therefore
\begin{equation}
 t_2-t_1\le \frac{4\pi c_\delta}{\delta\Qprod}M(u)^2 .
                                                               \label{eq:smoothbound}
\end{equation}
In particular, if $u\in C(I;\Sch(\R^3))$ satisfies \eqref{eq:hyp},
taking the supremum of \eqref{eq:smoothbound} over all compact
subintervals $[t_1,t_2]\subset I$ containing $0$ yields
\eqref{eq:lifespan} with $C_\delta=4\pi c_\delta/(\delta\Qprod)$.

We next extend this bound to general initial data $u_0\in H^1(\R^3)$
satisfying \eqref{eq:hyp}. Choose a sequence of Schwartz data
$u_{0,n}\to u_0$ in $H^1(\R^3)$, let $u_n$ be the corresponding maximal
solutions, and write $M_n=M(u_{0,n})$, $E_n=E(u_{0,n})$, and
$\delta_n=1-(M_nE_n)_+/(M(Q)E(Q))$. By the Sobolev embedding
$H^1(\R^3)\hookrightarrow L^4(\R^3)$, we have $M_n\to M(u)$,
$E_n\to E(u)$, and $\|\nabla u_{0,n}\|_{L^2_x}\to\|\nabla u_0\|_{L^2_x}$,
which implies $\delta_n\to\delta$. By the strict inequalities in
\eqref{eq:hyp}, for all sufficiently large $n$ we have
\[
 M_nE_n<M(Q)E(Q),\qquad
 \|u_{0,n}\|_{L^2_x}\|\nabla u_{0,n}\|_{L^2_x}>\Qprod,
\]
so each such $u_n$ satisfies \eqref{eq:hyp} with $\delta_n\in(0,1]$.

Fix a compact interval $[t_1,t_2]\subset I$ containing $0$, on which
$\|u(t)\|_{H^1_x}$ is bounded. Iterating the local $H^1$ well-posedness
and continuous dependence \cite{Caz} over finitely many subintervals
starting from $0$ shows that, for all sufficiently large $n$, the
solution $u_n$ exists on $[t_1,t_2]$ and satisfies $u_n\to u$ in
$C([t_1,t_2];H^1(\R^3))$.

To apply \eqref{eq:smoothbound}, we verify that
$u_n\in C([t_1,t_2];\Sch(\R^3))$ for each such $n$. By persistence of
Sobolev regularity for the cubic NLS \cite{Caz}, we have
$u_n\in C([t_1,t_2];H^2(\R^3))$, and hence
$u_n\in L^\infty([t_1,t_2];L^8(\R^3))$ by the Sobolev embedding
$H^2(\R^3)\hookrightarrow L^8(\R^3)$. For each integer $k\ge2$, define
the weighted space
\[
 \Sigma^k:=\bigl\{v\in H^k(\R^3):\ x^\alpha v\in L^2(\R^3)
 \text{ for all }|\alpha|\le k\bigr\}.
\]
Applying the weighted persistence theorem \cite[Theorem~4]{HNT} to the
cubic equation with $L^q$-exponent $q=8$ and equal derivative and weight
orders $k\ge2$, we obtain $u_n\in C([0,t_2];\Sigma^k)$ for all $k\ge2$. Since
$\bigcap_{k\ge2}\Sigma^k=\Sch(\R^3)$ with its usual topology, it follows
that $u_n\in C([0,t_2];\Sch(\R^3))$. By time reversal, the same holds on
$[t_1,0]$, and thus $u_n\in C([t_1,t_2];\Sch(\R^3))$.

Applying \eqref{eq:smoothbound} to $u_n$ on $[t_1,t_2]$, we obtain
\[
 t_2-t_1\le \frac{4\pi c_{\delta_n}}{\delta_n\Qprod}M_n^2 .
\]
Letting $n\to\infty$ and using $M_n\to M(u)$, $\delta_n\to\delta$, and
the continuity of $\delta\mapsto c_\delta$ on $(0,1]$, we obtain
\eqref{eq:smoothbound} for $u$ on $[t_1,t_2]$. Taking the supremum of
$t_2-t_1$ over all compact intervals $[t_1,t_2]\subset I$ containing
$0$ yields \eqref{eq:lifespan} with
$C_\delta=4\pi c_\delta/(\delta\Qprod)$. In particular, $T_*<\infty$
and $T^*<\infty$.

Finally, if $\|\nabla u(t)\|_{L^2_x}$ remained bounded along a sequence
of times tending to $T^*$ or $-T_*$, then by mass conservation
$\|u(t)\|_{H^1_x}$ would also be bounded along this sequence. By the
local $H^1$ theory, $u$ could then be extended beyond the endpoint,
contradicting the maximality of $I$. Hence
$\|\nabla u(t)\|_{L^2_x}\to\infty$ both as $t\uparrow T^*$ and as
$t\downarrow -T_*$.
\end{proof}

\section*{Acknowledgements}

The author would like to thank Professor Bin Liu and Lei Zhang for their constant encouragement, generous support, and many patient discussions throughout the preparation of this work.

\section*{Declarations}

\noindent\textbf{Funding.}
This work was partially supported by the National Key Research and Development Program of China (Grant No.~2023YFC2206100), and the Fundamental Research Funds for the Central Universities (Grant No.~2026BRSXB003).

\medskip
\noindent\textbf{Competing interests.}
The author has no competing interests to declare that are relevant to
the content of this article.

\medskip
\noindent\textbf{Data availability.}
Data sharing is not applicable to this article as no datasets were
generated or analyzed during the current study.


\begin{thebibliography}{99}

\bibitem{BF}
J.~Bellazzini and L.~Forcella,
\emph{Dynamical collapse of cylindrical symmetric dipolar Bose--Einstein
condensates},
Calc. Var. Partial Differential Equations 60 (2021), no.~6, Paper No.~229.

\bibitem{CFR}
L.~Campos, L.~G. Farah, and S.~Roudenko,
\emph{Threshold solutions for the nonlinear Schr\"odinger equation},
Rev. Mat. Iberoam. 38 (2022), no.~5, 1637--1708.

\bibitem{Caz}
T.~Cazenave,
\emph{Semilinear Schr\"odinger Equations},
Courant Lecture Notes in Mathematics, 10, New York University, Courant
Institute of Mathematical Sciences, New York; American Mathematical
Society, Providence, RI, 2003.

\bibitem{CW}
T.~Cazenave and F.~B. Weissler,
\emph{The Cauchy problem for the nonlinear Schr\"odinger equation in $H^1$},
Manuscripta Math. 61 (1988), no.~4, 477--494.

\bibitem{Coffman}
C.~V. Coffman,
\emph{Uniqueness of the ground state solution for $\Delta u-u+u^3=0$ and a
variational characterization of other solutions},
Arch. Rational Mech. Anal. 46 (1972), 81--95.

\bibitem{DWZ}
D.~Du, Y.~Wu, and K.~Zhang,
\emph{On blow-up criterion for the nonlinear Schr\"odinger equation},
Discrete Contin. Dyn. Syst. 36 (2016), no.~7, 3639--3650.

\bibitem{DHR}
T.~Duyckaerts, J.~Holmer, and S.~Roudenko,
\emph{Scattering for the non-radial 3D cubic nonlinear Schr\"odinger
equation},
Math. Res. Lett. 15 (2008), no.~6, 1233--1250.

\bibitem{DR10}
T.~Duyckaerts and S.~Roudenko,
\emph{Threshold solutions for the focusing 3D cubic Schr\"odinger equation},
Rev. Mat. Iberoam. 26 (2010), no.~1, 1--56.

\bibitem{DR15}
T.~Duyckaerts and S.~Roudenko,
\emph{Going beyond the threshold: scattering and blow-up in the focusing
NLS equation},
Comm. Math. Phys. 334 (2015), no.~3, 1573--1615.

\bibitem{Glassey}
R.~T. Glassey,
\emph{On the blowing up of solutions to the Cauchy problem for nonlinear
Schr\"odinger equations},
J. Math. Phys. 18 (1977), no.~9, 1794--1797.

\bibitem{GI}
S.~Gustafson and T.~Inui,
\emph{Blow-up or Grow-up for the threshold solutions to the nonlinear
Schr\"odinger equation},
Dyn. Partial Differ. Equ. 20 (2023), no.~3, 213--225.

\bibitem{HNT}
N.~Hayashi, K.~Nakamitsu, and M.~Tsutsumi,
\emph{Nonlinear Schr\"odinger equations in weighted Sobolev spaces},
Funkcial. Ekvac. 31 (1988), no.~3, 363--381.

\bibitem{HRrad}
J.~Holmer and S.~Roudenko,
\emph{A sharp condition for scattering of the radial 3D cubic nonlinear
Schr\"odinger equation},
Comm. Math. Phys. 282 (2008), no.~2, 435--467.

\bibitem{HR}
J.~Holmer and S.~Roudenko,
\emph{Divergence of infinite-variance nonradial solutions to the 3D NLS
equation},
Comm. Partial Differential Equations 35 (2010), no.~5, 878--905.

\bibitem{KRRT}
E.~A. Kuznetsov, J.~Juul Rasmussen, K.~Rypdal, and S.~K. Turitsyn,
\emph{Sharper criteria for the wave collapse},
Physica D 87 (1995), nos.~1--4, 273--284.

\bibitem{OT}
T.~Ogawa and Y.~Tsutsumi,
\emph{Blow-up of $H^1$ solution for the nonlinear Schr\"odinger equation},
J. Differential Equations 92 (1991), no.~2, 317--330.

\bibitem{Weinstein}
M.~I. Weinstein,
\emph{Nonlinear Schr\"odinger equations and sharp interpolation estimates},
Comm. Math. Phys. 87 (1983), no.~4, 567--576.

\end{thebibliography}
\end{document}